\documentclass[11pt]{article}

\usepackage{graphicx,epsfig}
\usepackage{amsmath}

\usepackage{color}

\newcommand{\wavelength}{\lambda_0}

\graphicspath{
{./}
{./figu/}
}

\def\d{{\partial}}
\begin{document}
\vskip 2cm
\begin{center}
{\Large  \bf  SIMULATION OF LASER PROPAGATION IN A PLASMA WITH A FREQUENCY WAVE EQUATION}
\vskip 1cm
{ \Large S.\ Desroziers$^{\dag}$,  F.Nataf $^{\ddag}$ , {R.\ Sentis}$^{\dag}$  }
\vskip 1cm
  {$^{\dag}$ CEA/Bruyeres, 91680 Bruyeres, France
    \\  remi.sentis@cea.fr,\ \  sylvain2.desroziers@cea.fr \\
$^{\ddag}$ Labo J-L-Lions, CNRS UMR 7598, Universit\'e Paris VI, 75013 Paris, France\\
nataf@ann.jussieu.fr    }
\end{center}

\pagestyle{empty}
\thispagestyle{empty}

\bibliographystyle{plain}

\bigskip

\begin{abstract} 
The aim of this work is to perform numerical simulations of the propagation of a laser in a plasma. At each time step, one has to solve a Helmholtz equation in a domain which consists in some hundreds of millions of cells. To solve this huge linear system, one uses a iterative Krylov  method with a preconditioning  by a separable matrix. The corresponding linear system is solved with a block cyclic reduction method. Some enlightments on the parallel implementation  are also given. Lastly, numerical results are presented including some features concerning the scalability of the numerical method on a parallel architecture. 
\end{abstract}

\section{Introduction}
\label{sec:intro}

The numerical simulation of propagation of high power intensity lasers in a
plasma is of importance for the \textquotedblright NIF
project\textquotedblright\ in USA and \textquotedblright LMJ Facility
project\textquotedblright\ in France. It is a very challenging area for
scientific computing indeed the laser wave length $2\pi /k_{0}$ is equal to a
fraction of one micron and the simulation domain has to be much larger than
500 microns. One knows that in a plasma the index of refraction is equal to $\sqrt{1-N(\mathbf{x})},$  where
$N(\mathbf{x} ) =N_{e}(\mathbf{x})/N_{c}$ and
$N_{e}$ is the electron plasma density at
position $\mathbf{x}$ and the critical density $N_{c}$ is a constant depending only on the wave length. Of course, the laser  propagates only in the region where  $N(\mathbf{x}) \le 1 .$
In macroscopic simulations (where the simulation lengths are in the order of
some millimeters), geometrical optics models are used and numerical
solutions are based on ray tracing methods. To take into account more specific
phenomena such as diffraction, autofocusing and filamentation, one generally
uses models based on a paraxial approximation of the full  Maxwell equations. This kind approximation is based on the assumption that  the  density $N(\mathbf{x})$ is close to a mean value $N_{av}$  ; it allows to make an expansion ok W.K.B. type  with a constant wave vector.  At the end of section 2, we recall the paraxial equation (\ref{parax}) ;
see for example \cite{dorr}, \cite{lesdi} in a classical framework and \cite{doum} for an analysis of this equation in a tilted frame and a numerical method.   But, there are situations where the macroscopic variations of the plasma density $
N_{e}$ are not small, particularly  in the zones which are just before the critical density. In this zone, the laser beam
 undergoes a strong change of direction near a surface called \emph{caustic surface}. That is to say the
wave vector is strongly varying near this surface, the paraxial
approximation is no more valid  and one has to deal with a model based on a
  frequency wave equation (obtained by time envelope of the
  solution of the full Maxwell equations). The model is
described in the section
2. For a derivation of the models and a physical exposition of the phenomena under interest, see e.g. \cite{sen} or \cite{KOL}.

 This paper is aiming at
describing the numerical methods for solving the
frequency wave   equation  and the coupling with the model for the plasma behavior.   At each time step, one has to find the  solution $\psi=\psi (\mathbf{x} ) $ of the following Helmholtz problem

\begin{equation}
k_{0}^{-2}\Delta \psi +\left(
(1-N)+i \mu \right) \psi =f  \label{hhh}
\end{equation}%
where $f$ is a given complex function and
$\mu $ a positive function related to the absorption
of the laser by the plasma.

 In this paper, only 2D problems are considered but the method may
be extended to 3D simulations. Let us set $\mathbf{x}=(x,y)$ the two spatial coordinates. The key
assumption is that the gradient of the density $N(\mathbf{x})$ is parallel to the $x$-axis, then we set 
\begin{equation}
N(x,y)=N_{0}(x)+\delta N(x,y)  \label{deco}
\end{equation}%
where $N_{0}$ depends on the $x$ variable only and $\delta N$ is small
compared to 1. This allows to deal with real physical situations as it is shown in the numerical applications below.

 The simulation domain is a rectangular box and a classical finite difference method is used for the spatial discretization ; for an accurate solution it is necessary to have a spatial step equal
to a fraction of the wave length. If $n_{x}$ and $n_{y}$ denote the number
of discretization points in each direction, it leads to solve a the linear
system with $n_{x}n_{y}$ degrees of freedom (which may be in the order of $10^{8}$
for a typical 2D spatial domain).

One chooses an iterative method of Krylov
type with a preconditioning by a matrix corresponding to the discretization of (\ref{hhh})  with $N$  replaced by $N_0 .$
This leads to solve a linear system corresponding  to
a separable tri-diagonal block matrix  (each block is
a   $n_x \times n_x$ matrix), then a block cyclic reduction
method  may be used ; see for instance 
\cite{meu}  and  \cite{rossi} for this kind of method. 
The crucial point for this method is to decompose the
unknown  onto the basis of the $n_{x}$ eigenvectors of the
main-diagonal block matrix.

The paper is organized as follows. After the statement of the model in section 2, we present the main difficulties for the numerical simulation of such problems in section 3. In section 4, the numerical scheme for the Helmholtz solver is presented, especially the method for solving the preconditioner with the block cyclic reduction method ; some enlightments on the parallel implementation and on the coupling with the hydrodynamics part are also given. Lastly, numerical results are presented including some features concerning the scalability of the numerical method.

\section{Statement of the model}

Our goal is to
perform simulations 
taking into account diffraction, refraction and auto-focusing phenomena 
and it is necessary to perform a coupling between the 
 the fluid dynamics system for the plasma behavior and the frequency wave equation for the laser propagation (notice that the Brillouin parametric instabilities which create laser backscattering
are not taken into account up to now).

The laser beam is characterized by an electromagnetic wave with a fixed
pulsation, so it may be modeled by the
time envelope $\Psi =\Psi (t,\mathbf{x})$ of this electric field.
It is a slowly time varying complex function. On the other hand, for
modeling the plasma behavior one introduces the non-dimension electron
density $N=N(t,\mathbf{x})$ and the plasma velocity $\mathbf{U}=\mathbf{U}(t,%
\mathbf{x})$.

\textbf{Modeling of the plasma. }   For the plasma, the simplest model is the 
following fluid model. Let us denote $P=P(N, T_e)$ a smooth function of the density $N$  and of the electron temperature
$T_e$ ( $T_e$  is a very smooth given function of the spatial variable $\mathbf{x}$). Then one has to solve the following
barotropic Euler system~: 
\begin{eqnarray}
\frac{\partial }{\partial t}N+\nabla (N\mathbf{U}) &=&0,  \label{moli} \\
\frac{\partial }{\partial t}(N\mathbf{U})+\nabla (N\mathbf{UU})+\nabla
(P(N, T_e)) &=&-N\gamma _{p}\nabla |\Psi |^{2}.  \label{molii}
\end{eqnarray}

The term $\gamma _{p}\nabla |\psi|^{2}$ corresponds to a ponderomotive force
due to a laser pressure (the coefficient $\gamma _{p}$ is a constant
depending only on the ion species).

\textbf{Modeling of the laser beam. }
Let us denote $\epsilon = k_0^{-1}.$ The laser
field $\Psi =\Psi (t,\mathbf{x})$ is a solution to the following frequency wave equation 
(which is of Schr\"{o}dinger type)

\begin{equation}
2i\frac{1}{c}\frac{\partial }{\partial t}\Psi+\epsilon \Delta
\Psi+\frac1{\epsilon}(1-N) \Psi+i\nu \Psi=0,  \label{cham}
\end{equation}
where the absorption coefficient $\nu $ is a  real coefficient related
 to the absorption of the laser intensity by the plasma
and  $c$  the light speed.

\textbf{Boundary conditions. }  The laser beam is assumed to enter into the domain at $x=0.$ Denote by $
\mathbf{e}_{b}$  the unit vector in the direction of the incoming beam.
Since the density $N$ depends mainly on the $x-$variable, we may denote by $%
N^{in}$ the mean value of the incoming density on the boundary and by $%
N^{out}$ the mean value of the density on the outgoing boundary . The
boundary condition at $x=0$ reads (with $\mathbf{n}$ the outwards normal to the boundary) 
\begin{equation}
(\epsilon \mathbf{n}.\nabla
+i\mathbf{K}.\mathbf{n})(\Psi -\alpha
^{in}e^{ik_{0}\mathbf{K}\mathbf{x}})=0.  \label{bce}
\end{equation}%
where  $\mathbf{K}=\mathbf{e}_{b}\sqrt{1-N^{in}}$, $\alpha ^{in}=\alpha ^{in}(y)$ is a smooth function which is, roughly
speaking, independent of the time. On the part of the boundary $x=x_{\mathrm{max}}$,
 there are two cases according to the value $N^{out}$ :

i) If $N^{out}>1,$ the wave do not propagate up to the boundary and the
boundary condition may read as $\partial \Psi /\partial x=0.$

ii) If $N^{out}\leq 1,$ one has to consider a transparent boundary
condition. Here we take the simplest one, that is to say 
$
\epsilon \mathbf{n}.\nabla +i\sqrt{1-N^{out}})\Psi =0. 
$

On the other hand, on the part of the boundary corresponding to $y=0$ and $%
y=y_{\mathrm{max}}$, it is crucial to have a good transparent boundary
condition, so we introduce perfectly matched layers (the P.M.L. of \cite{pml}%
). For the simple equation $-\Delta \psi -\omega ^{2}\psi =f,$\ this
technique amounts to replace in the neighborhood of the boundary, the
operator $\frac{{\partial }}{{\partial }y}$ by $\left( 1+\frac{\sigma }{%
\mbox{i}\omega }\right) ^{-1}\frac{{\partial }}{{\partial }y},$ where $%
\sigma $ is a damping function which is not zero only on two or three wave
lengths and which increases very fast up to the boundary.
Notice that the feature of this method is that it is necessary to modify the
discretization of the Laplace operator on a
small zone near the boundaries.

{\bf The paraxial equation.} For the sake of completeness,
we recall now the paraxial 
approximation equation which is valid only if the plasma density is a very smooth function, in such a way that we can take  $N_0 =N_{\rm av}$  where $N_{\rm av}$ in a constant. So we can define a mean wave
vector   $\mathbf{K}=\mathbf{e}_{b}\sqrt{1-N_{av}} $
and the laser beam is now characterized by the space and time
 envelope of the electric field $E=E(t,\mathbf{x})$, that is to say  
$$
 \Psi(t,\mathbf{x})=E(t,\mathbf{x})
 e^{i\mathbf{K}.\mathbf{x}/\epsilon }.
$$ 
The envelope $E$  is assumed to be slowly varying with the space variable 
and thus satisfies

\begin{equation}
i\left( \frac{2}{c}\frac{\partial E}{\partial t}+2\mathbf{K.}\nabla
E+\nu E\right) +\epsilon(\Delta _{\bot
}^{K}E)-\frac1{\epsilon} (N-N_{av})E=0.  \label{parax}
\end{equation}
where
$\Delta _{\bot}^{K}$ denotes the
 Laplace operator in the hyperplane transverse
 to $K.$
 It is necessary to supplement  equation
(\ref{parax}) with  a boundary condition 
  on the incoming boundary which is $E(0,.)=\alpha^{in}$ (and an initial
  condition). See \cite{dorr}, \cite{lesdi} , \cite{doum}.

\section{Difficulties}
\label{sec:difficulties}

The discretization and the solving of the above system of partial differential equations is very challenging 
since different space scales are to be considered. On the other hand, it leads to very large linear system to solve.

\subsection{Multiscale in space}
\label{sec:multiscale}
    
For solving (\ref{cham}),
the spatial mesh has to be very fine,  
  $h_{Helmholtz}\simeq \wavelength/10$ or less
each direction (recall that  $\wavelength =2\pi / k_0 $) ; this mesh is called in the sequel the Helmholtz grid.
 But the modulus $|\psi |$ of the electric field is
slowly varying with respect to the spatial variable,  thus one can use a coarse
mesh for the simulation of the Euler system,  typically one can set  $h_{fluid}\simeq  \wavelength/2$.

 For the numerical solution of the fluid system,
we refer to the method described in \cite{jour} or \cite{lesdi} which has been
implemented in a parallel platform called HERA,  the ponderomotive force is taken into account by adding the ponderomotive force  
 proportional to $ \nabla |\psi |^2$  to the  pressure force. The plasma density $N$, the velocity $\mathbf{U}$ and the laser intensity $  |\psi |^2$
are evaluated at the center of each cell. 

So we handle a two-level mesh of finite difference type : in a 2D
simulation, each cell of the fluid system is divided into $p_{0}\times p_{0}$
cells for the Helmholtz level, with $p_{0}=5$ or more.
At each time step $\delta t$ determined by the CFL criterion for the Euler
system,  one  has  to solve the frequency wave equation (\ref{cham}). For the time
discretization of this equation, an implicit scheme is used. The length $c\delta t$ is very large compared to the spatial step therefore the time
derivative term may be considered as a perturbation.
So, at each time step, if $\psi^{ini}$ denotes the value of the solution at the beginning of the time step, one has to find $\psi$ solution of the
following equation of the Helmholtz type 
\begin{equation}
 \epsilon^{2} \Delta \psi +\left((1-N) + i \mu _{1} \right) \psi =i
 \mu
_{0}\psi ^{ini}  \label{hee}
\end{equation}%
where $\mu _{0}= \epsilon   2/(c\delta t),$
 $\mu_1  =  \epsilon (  2/(c\,\delta t)  +\nu ).$ 
 This equation is supplemented by a boundary
 condition at $x=0$
 \begin{equation}
\epsilon \frac {\partial}{\partial
x}+i\mathbf{K}_x) \psi=
\epsilon \frac {\partial}{\partial x}+i\mathbf{K}_x) 
 (\alpha
^{in}e^{i\mathbf{K}_y  . y/\epsilon}). 
\label{bbb}
\end{equation}%

and another one in $x=x_{\mathrm{max}}$ as above.

\subsection{Large Scale Problem}
\label{sec:largescaleproblem}
As we shall see, the main difficulty comes from the Helmholtz equation: the number of unknowns is quite large and the properties of the resulting linear system makes 
it hard to solve. The linear system is symmetric but not Hermitian ; these properties are inherited by the discretized equations. 
The resulting linear systems are thus difficult to precondition. 
Concerning the preconditioning, the hypothesis (\ref{deco}) leads to replaced the original system by another one which is
simpler since it does not take into account the perturbation $\delta N(x,y).$ The corresponding linear system to be solved
leads to a five-diagonal symmetric non-hermitian matrix $A_{G}$ 
\begin{equation}
   \label{eq:cyclicmatrix}
A_{G}=\left( 
\begin{array}{cccccc}
\ \beta +A & -T &  &  &  &  \\ 
-T & A & -T &  &  &  \\ 
& -T & A & -T & ... &  \\ 
&  & -T & A & ... & 
\end{array}%
\right) 
\end{equation}%
where $T$ is equal to a constant times the identity
matrix of dimension $n_{x},$ the matrix $A$ of dimension 
$n_{x}$ and  $\beta$ is a complex constant.
The matrix $A_{G}$ is separable and therefore a block cyclic reduction method may be used for its numerical solution.

Notice that for a realistic simulation where the sizes of the domain is $500$ $\mu m$
 times $700$ $\mu m$ and the wave length equal to $0.35$ $\mu m,$ it leads to $12$ millions fluid unknowns and $300$ millions  unknowns on the Helmholtz grid .\\

\section{Numerical Strategies}
\label{sec:numstrat}

\subsection{Helmholtz solver}
\label{sec:solver}

Due to the size of the problem, a direct solver (even a parallel one) can not be used. We have to use an iterative method, in our case a preconditioned Krylov solver. 
As for the preconditioner, it seems difficult to propose one which would
be valid booth in the central zone where a pure Helmholtz equation is
considered and in the Perfectly Matched Layers (PML). Therefore, we first make a decomposition of the domain into three subdomains: two thin PML layers and a 
large central domain with appropriate interface conditions, see \S~\ref{sec:ddm}. The preconditioning of the central problem is presented in \S~\ref{sec:cyclic}, it corresponds to the solution to 
an approximated equation. 

 Let us mention right away that a multigrid method for the Helmholtz equation \eqref{eq:1:Fnat} would not work. Indeed, multigrid 
methods are efficient only if a large enough damping parameter is present in the
equation. Here, the coefficient $c \delta t $ is
larger than $100$  wave lengths
and the term $\nu$ is typically 
given by the formula 
\[
\nu = \nu_C N_0(x)^2
\]
with $1 / \nu_C $ in the order  of $15 \mu m  $  (typical values of
the density are
around $0.4$). So the damping term $ k_0 \mu_1 = k_0 ( \mu_0 + \nu ) $ is quite small 
when compared to the wave length  $2 \pi k_0 ^{-1}$ and it is too weak for a multigrid method.

   To solve the linear system arising from the discretization of the Helmholtz equation 
\begin{equation}\label{eq:1:Fnat}
 \epsilon^{2} \Delta \psi + \mbox{i} \mu_1 \psi +
      (1-N_{0}(x)) \psi - \delta_{N}(x,y) \psi  =  \mbox{i} \mu_0 \psi^{ini},
\end{equation}
we use the fact that the function $N_0$  depends only of a one-dimension variable
and we deal with a Krylov method with a preconditioner which corresponds to a  separable
matrix. This preconditioner may be interpreted as the discretization of  the operator 
\begin{equation}\label{eq:2:Fnat}
\psi \quad \mapsto  \quad \epsilon^{2} \Delta \psi + \mbox{i} \mu_0 \psi +
      (1-N_{0}(x)) \psi\,.
\end{equation}
with the same boundary conditions.
A block cyclic reduction method is then used for solving the corresponding linear system. 

Let us mention that the idea of preconditioning a variable coefficient Helmholtz equation by a problem amenable to the separation of variables technique was 
investigated in \cite{Plessix:2003:SVP} in the context of seismic modeling. Although in this context the method was not satisfactory, we shall see that in our context (laser-plasma interaction), results are indeed very good. We mention as well another method  \cite{Erlangga:2006:CMI} based on the preconditioning of \eqref{eq:1:Fnat} by a shifted Helmholtz equation with $\beta_0 \simeq 0.3-0.5$
\begin{equation*}
 \epsilon^{2} \Delta \psi + \mbox{i} \mu_1  \psi +
     (1+\mbox{i}\beta_0) [(1-N) \psi ]
\end{equation*}
that would be amenable to a multigrid solver. 

\subsubsection{Domain Decomposition}
\label{sec:ddm}
\begin{figure}
\centering
 \includegraphics[width=7cm]{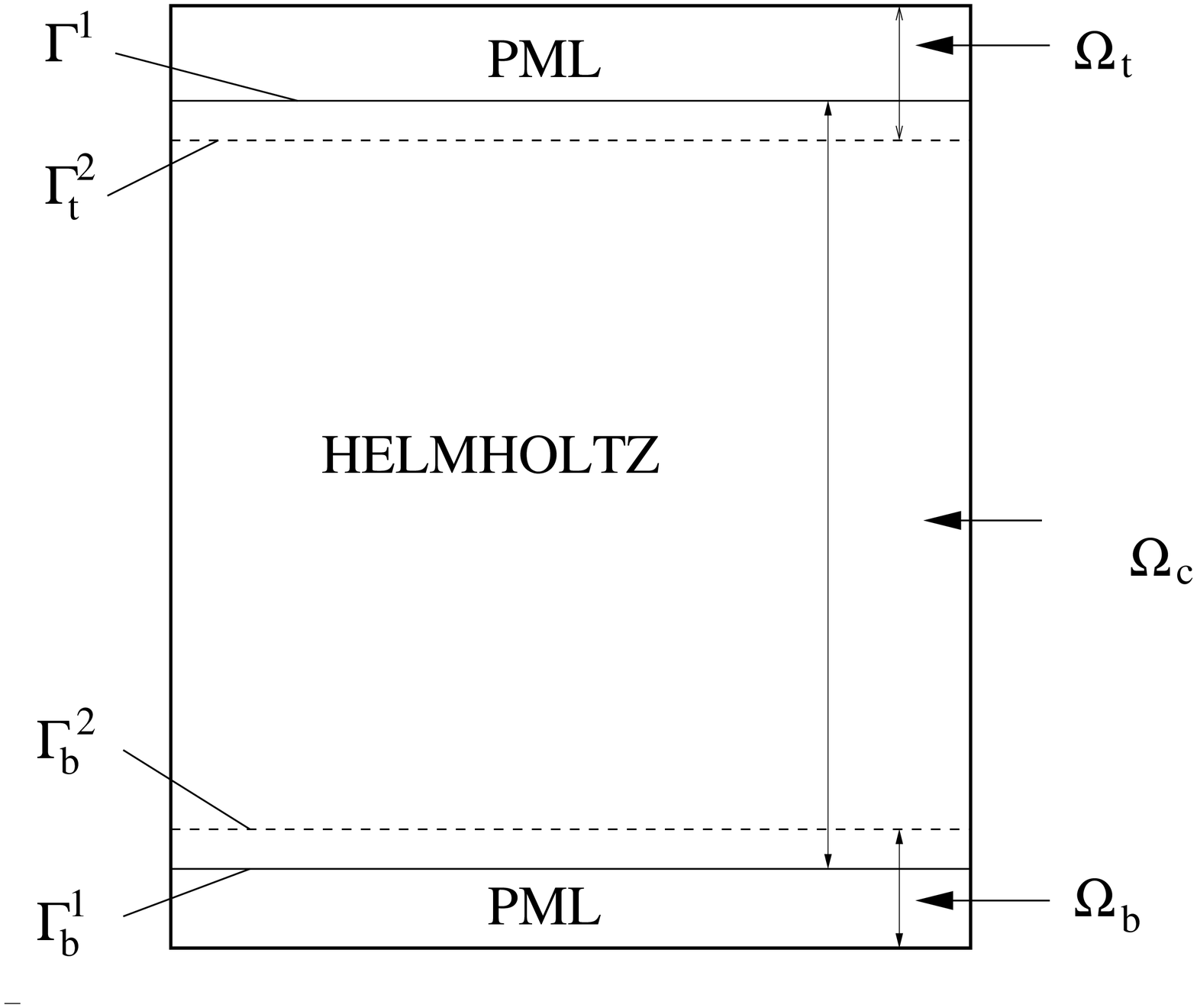}
\caption{Domain decomposition into three overlapping subdomains 
\label{fig:decomposition}}
  \end{figure}

The computational domain is divided into three overlapping subdomains: a purely Helmholtz zone $\Omega_c$ and two zones bordering 
it above $\Omega_t$ and below $\Omega_b$, see fig. \ref{fig:decomposition}. In $\Omega_t$ and $\Omega_b$, we have both a PML and a
 Helmholtz region. 

The coupling between the subdomains is made via Robin interface conditions,
see \cite{Lions:1990:SAM} and \cite{Benamou:1997:DDH}. So
 solving of equation (\ref{hee}) leads to consider 
 the following coupled system of equations, where the unknown functions are $\psi_t ,\psi_{c}, \psi_{b}$

   $$ 
    \left\{
    \begin{array}{lcl}
      \displaystyle \epsilon^{2} \left[ \eta(y) \frac{\d}{\d y} \left( \eta(y) \frac{\d}{\d y} \right)
        + \frac{\d^{2}}{\d x^{2}} \right] \psi_{t} + {\rm i} \mu_1 \psi_{t} +
      (1-N_{0}) \psi_{t} = 0 & \mbox{in} & \Omega_{t}\\
      \displaystyle \frac{\d \psi_{t}}{\d y} + {\rm i} \alpha \psi_{t} = \frac{\d \psi_{c}}{\d y} +{\rm i}  \alpha \psi_{c}
      & \mbox{on} & \Gamma_{h}^{2}
    \end{array}
    \right.
$$
$$
    \left\{
    \begin{array}{lcl}
      \displaystyle \epsilon^{2} \Delta \psi_{c} + {\rm i}  \mu_1 \psi_{c} +
      (1-N_{0}) \psi_{c} - \delta_{N} \psi_{c} = 
      \mbox{i} \mu_0 \psi^{ini} & \mbox{in} & \Omega_{c}\\
      \displaystyle \frac{\d \psi_{c}}{\d y} + {\rm i} \alpha \psi_{c} = \frac{\d \psi_{t}}{\d y} +{\rm i}  \alpha \psi_{t}
      & \mbox{on} & \Gamma_{t}^{1}\\
      \displaystyle -\frac{\d \psi_{c}}{\d y} + {\rm i} \alpha \psi_{c} = -\frac{\d \psi_{b}}{\d y} +{\rm i}  \alpha \psi_{b}
      & \mbox{on} & \Gamma_{b}^{1}
    \end{array}
    \right.
$$
$$
    \left\{
    \begin{array}{lcl}
      \displaystyle \epsilon^{2} \left[ \eta(y) \frac{\d}{\d y} \left( \eta(y) \frac{\d}{\d y} \right)
        + \frac{\d^{2}}{\d x^{2}} \right] \psi_{b} + \mbox{i} \mu_1  \psi_{b} +
      (1-N_{0}) \psi_{b} = 0 & \mbox{in} & \Omega_{b}\\
      \displaystyle -\frac{\d \psi_{b}}{\d y} + {\rm i} \alpha \psi_{b} = -\frac{\d \psi_{c}}{\d y} +{\rm i}  \alpha \psi_{c}
      & \mbox{on} & \Gamma_{b}^{2}
    \end{array}
    \right.
  $$

The coupling interface condition are of Robin type with  $\alpha=0.5 / \epsilon .$ 
  The boundary condition at $x=0$ is given (\ref{bbb}) for $\psi_c$ and by $(\epsilon \frac {\partial}{\partial
x}+ \mbox{i} \mathbf{K}_x) \psi=0 \,$  for $\psi_t$  and $\psi_b.$ 
These equations are discretized by a finite difference scheme. Let us denote by $\Psi_t$, $\Psi_c$ and $\Psi_b$ the corresponding unknown vector in domains $\Omega_t$, $\Omega_c$ and $\Omega_b.$ The linear system to be solved reads~:
  \begin{equation}\label{eq:linear_system}
 M \left( \begin{array}{c} \Psi_t\\ \Psi_c \\ \Psi_b \end{array} \right) = \left( \begin{array}{c} b_t\\ b_c \\ b_b \end{array} \right)  
\ \ \mbox{ where }\ \ 
M = \left[
      \begin{array}{ccc}
        A_{P1} & C_{1} & 0 \\
        C_{2} & A_{H} & C_{3}\\
        0 & C_{4} & A_{P2}
      \end{array}
      \right]\,.
  \end{equation}
The blocks $(C_i)_{1\le i\le 4}$ are related to the discrete Robin interface conditions. 

\subsubsection{The matrix system}
\label{sec:cyclic}
Denote by $A_{\delta N}$  the diagonal matrix
corresponding to the discretization of the operator of
multiplication by $-\delta_N + {\rm i} (\mu_1 - \mu_0 ) $ and by $A_G$ the one corresponding to the discretization of
equation \eqref{eq:2:Fnat}, that is to say $A_{H} =A_{ G} + A_{\delta N},$ then the matrix $M$ may be decomposed as
  $$
M=M_D + M_E, 
$$
with
$$
 \begin{array}{ccc}    M_{D} = \left[
      \begin{array}{ccc}
        A_{P1} & 0 & 0 \\
        0 & A_{G} & 0\\
        0 & 0 & A_{P2}
      \end{array}
      \right]
    & \mbox{and} &
    M_{E} = \left[
      \begin{array}{ccc}
        0 & C_{1} & 0 \\
        C_{2} & A_{\delta N} & C_{3}\\
        0 & C_{4} & 0
      \end{array}
      \right]
  \end{array}
  $$

The principle of solving the linear system \eqref{eq:linear_system} is to use a Krylov method preconditioned by $M_D$ which is a block diagonal matrix ; here we choose a GMRES algorithm (without restarting since the number of iterations is quite low, see below section 5.2). To apply the  preconditioner, since matrices $A_{P1},\quad A_{P2}$ are small ones, they can be factorized by a direct method and from the computational point of view, the key step is to use a fast solver for the matrix $A_G.$

 Let us describe the structure of the matrix $A_G$ coming from the discretization of equation \eqref{eq:2:Fnat}.  Denote by $\delta y$ is the space step in $y,$ by $\mathbf{I}$ the identity matrix of dimension $n_x .$ Let $A_0$ be the symmetric tridiagonal
  matrix which corresponds to the discretization of the following 1D problem 
  \begin{equation*}
 \epsilon^{2} 
\frac{\partial ^2}{\partial  x^2}  \psi  +
      (1-N_{0}(x)) \psi,
\end{equation*}
with the boundary condition \eqref{bbb}.
Their coefficients are  real except the one in the first line and the first
  column (due to the boundary condition). Then we set
 \begin{eqnarray*}
&&
A =  A_0 + i \mu_0  \mathbf{I} -  \frac{2 \epsilon ^2 }{ \delta y^{2} }\mathbf{I}
\\
&&
T= - \frac{\epsilon ^2 }{ \delta y^{2} } \mathbf{I}, \qquad
\\
&& B=A_0 + i (\mu_0 - \frac{\alpha \epsilon ^2 }{ \delta y } )  \mathbf{I} -  \frac{\epsilon ^2 }{ \delta y^{2} } \mathbf{I} = 
A + \beta \mathbf{I}
\end{eqnarray*}
where $\beta= - i {\alpha \epsilon ^2 }/{ \delta y }   +  {\epsilon ^2 }/{ \delta y^{2} } $ .
Now, if we denote  $ \Psi_c=(u_1, u_2,  ... u_{n_y})$ and $ f=(f_1, f_2,  ... f_{n_y})$ where the elements $u_m$ and $f_m$ are $n_x -$vectors, the system $A_G \Psi_c =f$ reads as follows
\begin{equation}\label{probleme_initial_RC}
    \left(
    \begin{array}{ccccc}
     B & -T & & & \\
      -T &    A   &   -T   &        &     \\
      & \ddots & \ddots & \ddots &        \\
      &        &   -T   &    A   &   -T   \\
      &        &        &   -T   &   B
    \end{array}
    \right)
    \left(
    \begin{array}{c}
      u_{1}   \\
      u_{2}   \\
      \vdots  \\
      u_{n_{y}-1} \\
      u_{n_{y}}
    \end{array}
    \right)
    =
    \left(
    \begin{array}{c}
      f_{1}   \\
      f_{2}   \\
      \vdots  \\
      f_{n_{y}-1} \\
      f_{n_{y}}
    \end{array}
    \right)
  \end{equation}

\subsubsection{Cyclic Reduction}

In order to solve system \eqref{probleme_initial_RC}, we use the block cyclic
reduction method.  Let us recall the principle of this method, assuming that $n_{y} = 2^{k} - 1$ for the
  sake of simplicity.  We know that  $A$ and
  $T$ are commutative. 
  Consider 3 successive lines of (\ref{probleme_initial_RC}) for
  $i=2,4,...,n_{y}-1$ :
  \begin{equation}\label{systeme}
    \left\{
    \begin{array}{ccccccccccc}
      - Tu_{i-2}& + & Au_{i-1} & - & T u_{i} &   &           &   &           & = & f_{i-1} \\
      & - & Tu_{i-1} & + & A u_{i} & - & T u_{i+1} &   &           & = & f_{i}   \\
      &   &          & - & Tu_{i}  & + & A u_{i+1} & - & T u_{i+2} & = & f_{i+1}.
    \end{array}
    \right.
  \end{equation}
After a linear combination of these lines, we get :
  \begin{equation}\label{reduit}
    - T^{2}A^{-1}u_{i-2} + \left(A - 2 T^{2}A^{-1}\right) u_{i} - T^{2}A^{-1}u_{i+2} = f_{i} +
    T A^{-1} \left(f_{i-1} + f_{i+1} \right)
  \end{equation}

After this first step, the elimination procedure may be performed
again by induction. That is to say,  denote $A^{(0)} = A$, $B^{(0)} =
B$, $T^{(0)} = T$ et $f^{(0)} = f\, ;$
after $r$ elimination steps, the reduced system for $0 \leq r \leq
k-1$ owns $2^{k-r} - 1$ blocs and reads as:
  $$
    \left(
    \begin{array}{ccccc}
      B^{(r)}  & -T^{(r)} & & &                \\
      -T^{(r)} &  A^{(r)} & -T^{(r)} &         \\
      & \ddots & \ddots   & \ddots   &         \\
      &        & -T^{(r)} & A^{(r)}  & -T^{(r)}\\
      &        &          & -T^{(r)} & B^{(r)} \\
    \end{array}
    \right)
    \left(
    \begin{array}{c}
      u_{2^{r}}   \\
      u_{2.2^{r}} \\
      \vdots      \\
      u_{(n_{y}-1) - 2^{r} + 1} \\
      u_{n_{y} - 2^{r} + 1}
    \end{array}
    \right) 
    =
    \left(
    \begin{array}{c}
      f^{(r)}_{2^{r}}   \\
      f^{(r)}_{2.2^{r}} \\
      \vdots  \\
      f^{(r)}_{(n_{y}-1) - 2^{r} + 1} \\
      f^{(r)}_{n_{y} - 2^{r} + 1}
    \end{array}
    \right)
  $$
  where for $r=1,...,k-2$ :
  \begin{eqnarray}
  \nonumber
    A^{(r)} & = & A^{(r-1)} - 2
    \left(T^{(r-1)}\right)^{2}\left(A^{(r-1)}\right)^{-1}\\  
    B^{(r)} & = & A^{(r-1)} -
    \left(T^{(r-1)}\right)^{2}  \left(   \left(A^{(r-1)}\right)^{-1} +
      \left(B^{(r-1)}\right)^{-1}  \right)
      \label{B} \\
      \nonumber
    T^{(r)} & = &
    \left(T^{(r-1)}\right)^{2}\left(A^{(r-1)}\right)^{-1}
  \end{eqnarray}

 For the right hand side, we get the induction formula  :

  \begin{equation}
  \label{sec}
    f^{(r)}_{i.2^{r}} = f_{i.2^{r}}^{(r-1)} +  T^{(r-1)}\left(A^{(r-1)}\right)^{-1}
      \left( f^{(r-1)}_{i.2^{r} - 2^{r-1}} + f^{(r-1)}_{i.2^{r} + 2^{r-1}} \right)
      \end{equation}
   
  After all the elimination steps, it remains only one equation
  for finding  $u_{2^{k-1}}.$ Once this value is obtained,  one
  deduces  all the other values  step by step recursively.
%  The method is illustrated by the sketch  on figure 2.

 % \begin{figure}[h]
 %%   \centering \epsfig{file=Images/reduction.eps,height=7cm,width=6cm}
  %  \caption{Steps of the Cyclic Reduction method  for 15 lines}
  %\end{figure}
  
  From a practical point of view, one has to perform these computations in the
  spectral basis of the eigenvectors of the matrix $A,$ which are of course also
  eigenvectors of $T,$ $B,$ $A^{(r)},$ $T^{(r)},$ $B^{(r)}$  for all $r.$
  
  \subsection{Parallel Implementation}
The implementation of the method has been made in the HERA platform, see \cite{lesdi},\cite{jour}.
  For the Helmholtz solver, one first have to find a orthonormal basis of  eigenvectors of $A_0.$  As a matter of
  fact,  even if the  matrix is not exactly real, we search a set of  eigenvectors
  which are orthogonal for the pseudo scalar product 
    $<u,v>=u^{T}\cdot v.$ To do this, we
   use the ``new  QD'' algorithm of Parlett (cf. \cite{Parlett:1995:NQD}) although 
    it was designed for Hermitian matrices.
   We conjecture that it is always possible to find such a basis of eigenvectors for our class of matrices. The only difficulty would be to find a non zero eigenvector $v$ such that  $v^{T}\cdot v =0,$ 
but in practice, we never encountered any problem  by using the  method.  In the method proposed in \cite{Parlett:1995:NQD}, which follows an idea of \cite{Rutishauser:1958:SEP}, the computation
  of the eigenvalues is based on a series of LU factorization of tridiagonal matrices. 
     This step is sequential but very cheap in terms of memory and CPU time requirements especially when compared to the QR algorithm which would manipulate full matrices. 
     Once the eigenvalues have been computed, the computation of the eigenvectors consists in finding the kernel of tridiagonal matrices. This task is distributed among the processors and is thus parallel. In our tests, this method was 40 times faster than the QR method.\\

So let us denote $Q$ the  matrix whose columns are the eigenvectors of $A^0 .$ The matrix $Q$ is orthonormal   for the pseudo scalar product, that is to say
  $$
  QQ^{T}=Q^{T}Q= {\mathbf I }
  $$
  
  Since  $T$ is the identity matrix up to a multiplicative constant,
  one can introduce the diagonal matrices $\Lambda^{(0)} $  and $\Gamma^{(0)} $
  \begin{equation}
    A = Q \Lambda^{(0)} Q^{T}, \ \ T = Q \Gamma^{(0)} Q^{T}.
  \end{equation}
  
  So we get 
  \begin{equation}\label{vapATr}
    A^{(r)} = Q \Lambda^{(r)} Q^{T}, \ \ T^{(r)} = Q \Gamma^{(r)} Q^{T}
  \end{equation}
with the following induction  formulas
  \begin{equation}
    \Lambda^{(r)} = \Lambda^{(r-1)} - 2 \left(
    \Gamma^{(r-1)}\right)^{2}\left(\Lambda^{(r-1)}\right)^{-1}, \quad
    \Gamma^{(r)} = \left( \Gamma^{(r-1)}\right)^{2}\left(\Lambda^{(r-1)}\right)^{-1}
    \label{vapTr}
  \end{equation}

 Let us summarize the algorithm

  $\bullet$ Introduce the vectors $\tilde{f}_{i}$ transformed of
  $f_{i}$ in the eigenvector
  basis 
    $$
    \begin{array}{lcr}
      \tilde{f}_{i} = Q^{T} f_{i} & \mbox{ for } & i = 1,\dots,n_{y}.
    \end{array}
    $$
    
  $\bullet$ At each step $r,$ the vector 
  $\tilde{f^r}_{i}$ transformed of ${f^r}_{i}$  of the  right hand
  side,  reads
    $$
    \tilde{f}^{(r)}_{i.2^{r}} = \tilde{f}_{i.2^{r}}^{(r-1)} +
    \Gamma^{(r-1)}\left(\Lambda^{(r-1)}\right)^{-1}
    \left( \tilde{f}^{(r-1)}_{i.2^{r} - 2^{r-1}} + \tilde{f}^{(r-1)}_{i.2^{r} + 2^{r-1}}
    \right)
    $$

  $\bullet$ One computes the vectors 
  $\tilde{u}_{2^{k-1}}$  by solving
    $$
    \Lambda^{(k-1)}\tilde{u}_{2^{k-1}} = \tilde{f}^{(k-1)}_{2^{k-1}}
    $$

  $\bullet$ One recursively  distributes the solutions by
  solving  sub-systems of the following type
    $$
    \Lambda^{(r)}\tilde{u}_{j.2^{r+1} - 2^{r}} = \tilde{g}^{(r)}_{j.2^{r+1} - 2^{r}}
    $$
    where
    $$
    \tilde{g}^{(r)}_{j.2^{r+1} - 2^{r}} =  \tilde{f}^{(r)}_{j.2^{r+1} - 2^{r}} + \Gamma^{(r)}
    \left( \tilde{u}_{(j-1).2^{r+1}} + \tilde{u}_{(j).2^{r+1}} \right)
    $$
    
  $\bullet$ Lastly, the  solution $u$ is given by
    $$
    \begin{array}{lcr}
      u_{i} = Q \tilde{u}_{i} & \mbox{ for } & i = 1,\dots,n_{y}.
    \end{array}
    $$

For the parallel implementation, the processors are shared out according to
horizontal slabs in a balanced way. Let us note that the first and last processors
contain PML layers and some lines of the central grid. The crucial point of the algorithm is the
multiplication of a full matrix $Q$ (and its transpose) of dimension $n_x \times n_x $ by
the set of $n_y $ vectors, i.e. a matrix-matrix multiplication. Within a realistic
framework for Symmetric MultiProcessors (SMP) architecture, the matrix $Q$ of a size of
several giga octets can be contained only in the memory of the nodes of
processors. This involves us to use a technique of hybrid parallelization of type
MPI-multithreading. One wishes to profit from the locality of the data between
processors of the same node and to use MPI for the others communications between
nodes. Cutting is carried out in order to avoid conflict of writing of the
threads.  %%%%EST CE CORRECT DU POINT DE VUE ANGLAIS ??? %%
For example, shearing in $N_{q} \times N_{b}$ threads leads to

  $$
  \left[
    \begin{array}{ccc}
      & Q_{1} & \\ 
      & \vdots & \\ 
      & Q_{N_{q}} & 
    \end{array}
    \right]
  \times
  \left[
    \begin{array}{ccc}
      & & \\
      B_{1} & \dots & B_{N_{b}} \\
      & & 
    \end{array}
    \right]
  = \left[
    \begin{array}{ccc}
      Q_{1} B_{1} & \dots & Q_{1} B_{N_{b}} \\ 
      \vdots & \ddots & \vdots \\ 
      Q_{N_{q}} B_{1} & \dots & Q_{N_{q}} B_{N_{b}} 
    \end{array}
    \right].
  $$
Each product $Q_{i} B_{j}$ is carried out by a call with the BLAS library. Let us
specify that the steps of computation in the spectral basis ({\it  i.e.} $\Lambda^{(r)} ,$ $\Gamma^{(r)} ,$  $ \tilde{f}^{(r)}_{i.2^{r}} ....$) are of a very negligible
cost in comparison of matrix-matrix products. Our tests show a good scalability of
the cyclic reduction and a great speed of execution.

\subsection{Coupling with the hydrodynamics part}
\label{sec:nonmatching}
Since the laser intensity $|\psi|^2$ is slowly varying according to the space
variable, one can deal with the ponderomotive force and the  hydrodynamic equations  on the coarse grid  (whose
size is  equal to one half of the wave length). To evaluate the laser intensity
$|\psi|^2$ on this coarse grid, it suffices to take the mean value of $|\psi|^2$
on the fine grid.
To get the electron density, i.e. $N_0$ and $\delta_N ,$ we use a linear
interpolation between the coarse grid and the fine grid.
Of course in the PML zone, one do not evaluate the ponderomotive force.

\section{Numerical Results}
\label{sec:numres}
\subsection{Test cases}

The boundary value $\alpha^{in}$ has to mimic laser beam ; to be realistic the
 profile of $\alpha^{in}$ corresponds to a juxtaposition of a lot of small hot spots , called {\it speckles} 
whose intensity is very high compared to the mean intensity of the beam.
The shape of a speckle is generally a Gaussian function whose width is about
a few micrometers.

One considers here a simulation domain of $100 \times 300$ wave lengths ; the initial profile of density is a linear function increasing from $0.1$ at $x=0$ to $1.1$ at $x=x_{max}$. The  profile of $\alpha^{in}$ contains only three speckles
At the Helmholtz 
level, one handles only 3 millions of cells. With 32 PEs (of the type EV67 HP-Compaq), the CPU time is only a few minutes per time step with approximatively 
10 Krylov iterations  at each time step.

\begin{figure}
\begin{center}
 \includegraphics[width=0.4\textwidth]{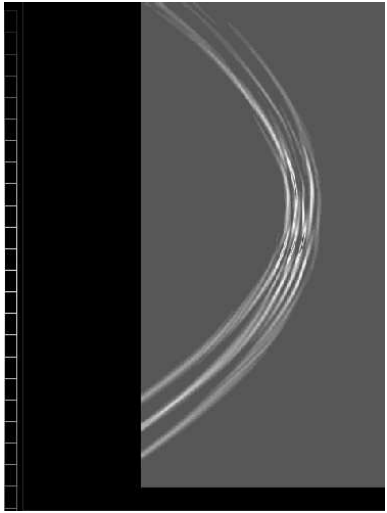}
  \caption{Laser intensity of a schematic  beamlet \label{DEU}}
  \end{center}
\end{figure}

 Without the coupling with the plasma, it is well known that the solution is very close to the one given by  Geometrical Optics  and corresponds to parallel speckles or beamlets which are curved and tangent to a 
 caustic line (this  line corresponds to  $x=x_{\star}$ such that $N_0(x_{\star} )=\cos ^2(\theta)$, where $\theta$ is the incidence angle of 
 the beamlets where they enter into the simulation box).
 With our model, if the laser intensity is small (which corresponds to a weak coupling
 with the plasma), one notices that a small digging of the plasma density occurs. 
 This digging is more significant when the 
 laser intensity is larger, then an autofocusing phenomenon takes place. On
 fig. \ref{DEU}, one sees the  map of the laser intensity that is to say the
 quantity $|\psi|^2$, which corresponds to this situation after some  time steps.
 We notice here that the three beamlets undergo autofocusing phenomena and something like a filamentation may be observed.

%\vspace{0.4cm}

\begin{figure}
\begin{center}
  \includegraphics[width=0.8\textwidth]{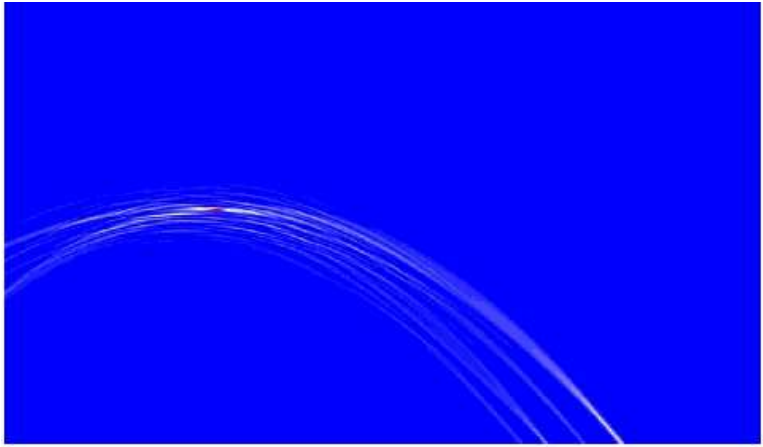}
  \caption{Laser intensity of a beamlet in a plasma after 6 ps
  without absorption. Notice the autofocusing of the beamlet near the
  caustic. The $x$ axis is here vertical. \label{TRO}}
  \end{center}
\end{figure}

Another case will be presented, corresponding to simulation domain of $700 \times 1200$ wave lengths. At the Helmholtz level, one handles 84 millions of 
cells and the simulation have run on 
128 PEs. The map of the laser intensity is shown on fig. \ref{TRO} after 6 picosecond
(corresponding  to about 15 time steps). We have set to zero the absorption
coefficient in order to have a  sharp problem.  The caustic lines
corresponds to about $N_0(x_{\star})=0.5.$ Here  the digging of the plasma is
locally very  important since the variation of density $\delta N$ 
reaches $0.07$  in a region where  $N_0(x)=0.45.$

\subsection{Numerical performance}
  \label{sec:scalabilityres}
\begin{figure}
  \label{fig:itcountvstime}
\centering
\includegraphics[height=7cm]{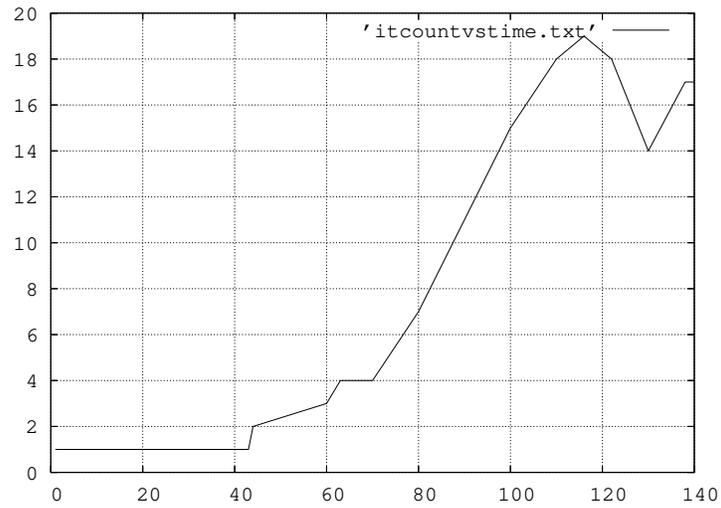}
\caption{Number of GMRES iterations versus time. \label{QUA}}
\end{figure}

We focus on the solving of the very large system \eqref{eq:linear_system} arising from the discretization of the Helmholtz equation by the preconditioned GMRES method presented above. In fig.~\ref{QUA}, we plot the iteration counts as a function of the physical time. As time increases, the density fluctuation $\delta _N$ gets larger and it is necessary to perform more iterations of the GMRES method. Nevertheless, we don't have more than 20 iterations. 

Let us address now the scalability of our computation method.
So we consider a fixed size problem of 40 million unknowns and we increase the number of processors. Table~\ref{tab:fixedsize} gives the elapsed time of one GMRES iteration. 
Due to the good parallel properties of the cyclic reduction, we see that the speed up is almost perfect.
When now the number of unknowns increases, one can easily check that the computational effort grows proportionally to $n_x^2 n_y$. In Table~\ref{tab:doubledsize}, from one column to the next, the number of points is doubled
 in each direction so the CPU time is 8 times larger.  Since the number of processors is multiplied by four, we check that the elapsed time is about two times larger. These two tables show that the parallelization of the cyclic reduction method work very well.

\begin{table}
\begin{center}
\begin{tabular}{||c||c|c|c|c||}
\hline 
    \mbox{Nb Procs} & 16 & 32 & 64 & 128 \\ \hline
    \mbox{elapsed time per GMRES it.} & 492s & 249s & 126s & 64s \\ \hline
    \mbox{Efficiency GMRES algo.} & 1 & 0.987 & 0.976 & 0.96 \\ \hline
\end{tabular} 
\caption{Scalability for a fixed size problem. \label{tab:fixedsize}}
\end{center} 
\end{table}

\begin{table}
\begin{center}
\begin{tabular}{||c||c|c|c|c|c||}\hline
    \mbox{Nb Procs} & 1 & 4 & 16 & 64 & 256 \\ \hline
    \mbox{\# d.o.f. $\times\,10^6$} & 0.4  & 1.6 & 6.3 & 25.4 & 101.6\\ \hline
    \mbox{elapsed time for QD algo.} & 1s & 3s & 12s & 48s & 189s \\ \hline
    \mbox{elapsed time per GMRES it.} & 4.8s & 11.6s & 24s & 47s & 93s\\ \hline
\end{tabular} 
\caption{Scalability for problems of increasing sizes. \label{tab:doubledsize}}
\end{center} 
\end{table}

\subsection{A more realistic case}

  In the realistic configurations, it may be useful to solve the paraxial equations on
  a part of the simulation domain in which this approximation is valid
  and the Helmholtz equations on the remainder where the validity of
  this approximation is no more true.  The mesh size for to numerical
  solution of  the paraxial equation is the same as for the fluid system (see \cite{jour}, \cite{dorr}). 
The coupling between the paraxial and Helmholtz parts is performed according to the classical boundary condition for the Helmholtz equation (\ref{bbb}) with $\alpha
^{in}$ replaced by $E^{out},$ which is the value of the solution to the paraxial equation at the interface boundary.
The advantage of the paraxial equation is that it can be solved by a marching method in space where only 1D systems have to be solved at each vertical line of unknowns ; see \cite{dorr}, \cite{doum}.

\begin{figure}
\hspace{4.5cm}
\centering
\includegraphics[width=9cm]{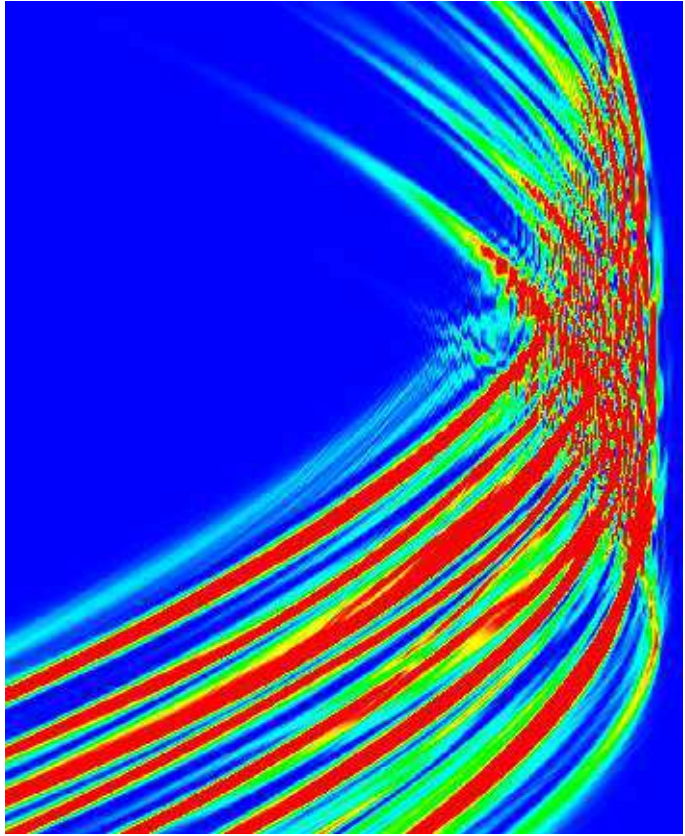}
  \caption{Zoom near the caustic line of the Laser intensity for a multi-speckle beam. \label{CIN}}
\end{figure}

Here we have performed a simulation with an initial density which is equal to $0.15$ up the third of the simulation domain and which ranges linearly up to $1$ at $x=x_{max}, $  the boundary value $\alpha^{in}$ mimics a multispeckle laser beam. We use the paraxial model in the third of the simulation domain and the frequency wave equation in the complementary part, then we have much less  
 unknowns to deal with for the Helmholtz problem. In this simulation, the computational domain size was $2000 \times
2000$ wave lengths. There were 200 million cells in the Helmholtz zone and   $4$    millions unknowns for the paraxial zone
The hydrodynamics equations are solved on a domain which consists in 16 million cells. 
 The computation was performed using 256
processors. The physical time of the simulation is equal to 11 $ps$. The elapsed time
for the full simulation was 8 hours.  On fig. \ref{CIN},  the map of the laser intensity is
represented at the end of simulation.

\begin{figure}
\hspace{4.5cm}
\centering
\includegraphics[width=9cm]{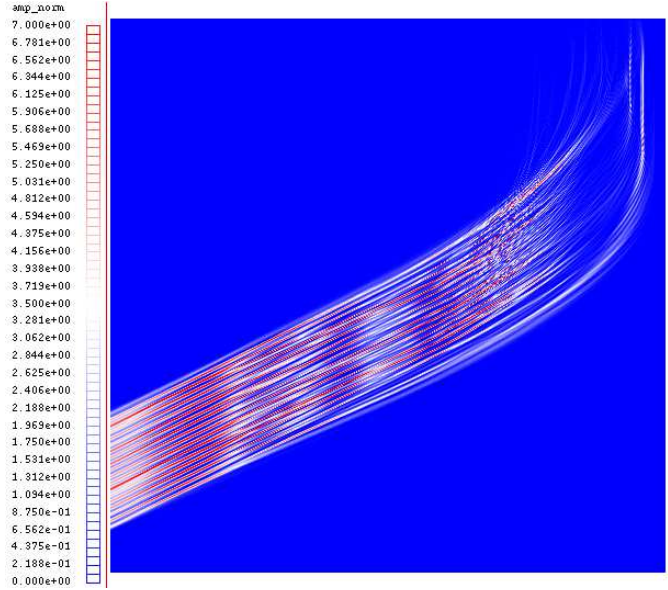}
  \caption{Laser intensity for a multiple beams laser for a plasma density ranging up to 1. \label{CIN2}}
\end{figure}

Lastly, we show on fig.~\ref{CIN2}, a zoom of the laser intensity near the caustic line in a numerical simulation of the same type than the previous one, except that the absorption coefficient is small but not zero. One can notice the great accuracy of simulation which shows interference patterns of the speckles.

\section{Conclusion and Prospects}
\label{sec:conclusion}
 In the framework of the hydrodynamics parallel platform  HERA,  we have developed a solver for the laser propagation based on the Helmholtz equation that can handle realistic computations on very large computational domains. The Helmholtz zone is coupled with a paraxial zone and a fluid plasma model. The assumption that the initial density $N$ depends mainly on the $x-$variable only allows to perform a preconditioning by a domain decomposition method (two PMLs and a large Helmholtz zone) where the linear system corresponding to the Helmholtz zone with a matrix $A_G$ may be solved efficiently by the block cyclic reduction method.  Most of the computer time is spent by applying a dense $n_x\times n_x$ matrix $Q$ to a set of $n_y$ vectors. We can thus achieve a very good scalability w.r.t to the number of unknowns in the $y$ direction. According to an increase of the fluctuation of the density, we notice an increase in the number of GMRES iterations per time step as the physical time increases, but this number remains small enough to have an acceptable CPU time.

In the future some CPU time may be saved by using a domain decomposition in the large central Helmholtz zone. For instance, simply dividing the Helmholtz zone in two vertical subdomains would decrease the size of the matrices $Q$ by a factor 4. Moreover the use of local (and thus more accurate) averages for the density in the preconditioner could break the increase in the number of iterations as the time increases. Another interesting strategy could be to not consider inside the inner iteration loop  of the Krylov method all the spatial domain that is to say all the $n_y$ vectors  but only the vectors which does not belong to some subinterval $[n_y^1, n_y^2]$ for instance the ones where the solution  varies very few  from an iteration to the other.

\bibliography{paperSDN}

\begin{thebibliography}{10}

\bibitem{lesdi}
Ph. Ballereau, M.~Casanova, Duboc F., Dureau D., Jourdren H., Loiseau P.,
  Metral J., O.~Morice, and Sentis R.
\newblock Simulation of the paraxial laser propagation coupled with
  hydrodynamics in 3{D} geometry.
\newblock {\em J. Scient. Comp. to appear}, 2007.

\bibitem{Benamou:1997:DDH}
J.~D. Benamou and B.~Despr\'{e}s.
\newblock A domain decomposition method for the {H}elmholtz equation and
  related optimal control.
\newblock {\em J. Comp. Phys.}, 136:68--82, 1997.

\bibitem{pml}
Jean-Pierre Berenger.
\newblock A perfectly matched layer for the absorption of electromagnetic
  waves.
\newblock {\em J. Comput. Phys.}, 114(2):185--200, 1994.

\bibitem{dorr}
R.M. Dorr, X.~Garaizar, and J.A.F Hittinger.
\newblock Simulation of laser-plasma filamentation using adaptive mesh
  refinement.
\newblock {\em J. Comput. Phys.}, 177:233--263, 2002.

\bibitem{doum}
Marie Doumic, Fran{\c{c}}ois Golse, and R{\'e}mi Sentis.
\newblock Un mod\`ele paraxial de propagation de la lumi\`ere: probl\`eme aux
  limites pour l'\'equation d'advection {S}chr\"odinger en coordonn\'ees
  obliques.
\newblock {\em C. R. Math. Acad. Sci. Paris}, 336(1):23--28, 2003.

\bibitem{Erlangga:2006:CMI}
Y.~A. Erlangga, C.~Vuik, and C.~W. Oosterlee.
\newblock Comparison of multigrid and incomplete {LU} shifted-{L}aplace
  preconditioners for the inhomogeneous {H}elmholtz equation.
\newblock {\em Appl. Numer. Math.}, 56(5):648--666, 2006.

\bibitem{KOL}
S.~H\"{u}ller, Mounaix Ph., D.~Pesme, and V.T. Tikhonchuk.
\newblock Interaction of two neighboring laser beams.
\newblock {\em Phys. Plasmas}, 4:2670--2680, 1997.

\bibitem{jour}
H.~Jourdren.
\newblock Hera hydrodynamics {AMR} plateform for multiphysics simulation.
\newblock In {\em {AMR} methods, theory and Applications, Plewa T, Linde T.,
  eds}. Lect. Notes Comp. Sciences, Springer, Berlin, 2005.

\bibitem{Lions:1990:SAM}
Pierre-Louis Lions.
\newblock On the {S}chwarz alternating method. {III:} a variant for
  nonoverlapping subdomains.
\newblock In Tony~F. Chan, Roland Glowinski, Jacques P{\'e}riaux, and Olof
  Widlund, editors, {\em Third International Symposium on Domain Decomposition
  Methods for Partial Differential Equations , held in Houston, Texas, March
  20-22, 1989}, Philadelphia, PA, 1990. SIAM.

\bibitem{meu}
G{\'e}rard Meurant.
\newblock A review on the inverse of symmetric tridiagonal and block
  tridiagonal matrices.
\newblock {\em SIAM J. Matrix Anal. Appl.}, 13(3):707--728, 1992.

\bibitem{Parlett:1995:NQD}
Beresford~N. Parlett.
\newblock The new qd algorithms.
\newblock In {\em Acta numerica, 1995}, Acta Numer., pages 459--491. Cambridge
  Univ. Press, Cambridge, 1995.

\bibitem{Plessix:2003:SVP}
R.~E. Plessix and W.~A. Mulder.
\newblock Separation-of-variables as a preconditioner for an iterative
  {H}elmholtz solver.
\newblock {\em Appl. Numer. Math.}, 44(3):385--400, 2003.

\bibitem{rossi}
Tuomo Rossi and Jari Toivanen.
\newblock A parallel fast direct solver for block tridiagonal systems with
  separable matrices of arbitrary dimension.
\newblock {\em SIAM J. Sci. Comput.}, 20(5):1778--1796 (electronic), 1999.

\bibitem{Rutishauser:1958:SEP}
Heinz Rutishauser.
\newblock Solution of eigenvalue problems with the {$LR$}-transformation.
\newblock {\em Nat. Bur. Standards Appl. Math. Ser.}, 1958(49):47--81, 1958.

\bibitem{sen}
R{\'e}mi Sentis.
\newblock Mathematical models for laser-plasma interaction.
\newblock {\em M2AN Math. Model. Numer. Anal.}, 39:275--318, 2005.

\end{thebibliography}

%%\begin{thebibliography}{99}
%%\end{thebibliography}

\end{document}